\documentstyle[11pt,epsf]{article}

\newlength{\mytopmargin}
\newlength{\myleftmargin}
\def\zz{\relax\hbox{\small \sf Z\kern-.4em Z}}

\newtheorem{lemma}{Lemma}
\newtheorem{prop}[lemma]{Proposition}
\begin{document}

\vspace{1cm}
\noindent
\begin{center}{\large \bf Random walks and random permutations}
\end{center}
\vspace{5mm}
\begin{center}
P.J.~Forrester
\end{center}

\vspace{.2cm}

\noindent
Department of Mathematics and Statistics, University of
Melbourne, Parkville 3052, Australia

\vspace{.2cm}

{\small
\begin{quote}
A connection is made between the random turns model of vicious walkers
and random permutations indexed by their increasing subsequences. Consequently
the scaled distribution of the maximum displacements in a particular
asymmeteric version of the model can be determined to be the same as
the scaled distribution of the eigenvalues at the soft edge of the GUE.
The scaling of the distribution gives the maximum mean displacement
$\mu$ after $t$ time steps as $\mu = (2t)^{1/2}$ with standard deviation
proportional to $\mu^{1/3}$. The exponent $1/3$ is typical of 
a large class of two-dimensional growth problems.
\end{quote}
}

Non-intersecting (vicious) random walkers were introduced into statistical
mechanics \cite{Fi84} as models of domain walls
and wetting in two-dimensional lattice systems, and have also received 
attention as exactly solvable systems \cite{Fo89e,Fo90b,Fo91f,BEO99,BO99}.
They can be viewed as directed lattice paths which start at sites say on
the $x$-axis and finish on sites on the line $y=n$, with the additional
constraint that the paths do not touch or overlap. Alternatively, vicious
 random walkers can be described as the stochastic evolution of particles
on a one-dimensional lattice, which at each tick of the clock move
to the left or to the right with a certain probability, subject to the
constraint that no two particles can occupy the one lattice site.

\begin{figure}
\epsfxsize=10cm
\centerline{\epsfbox{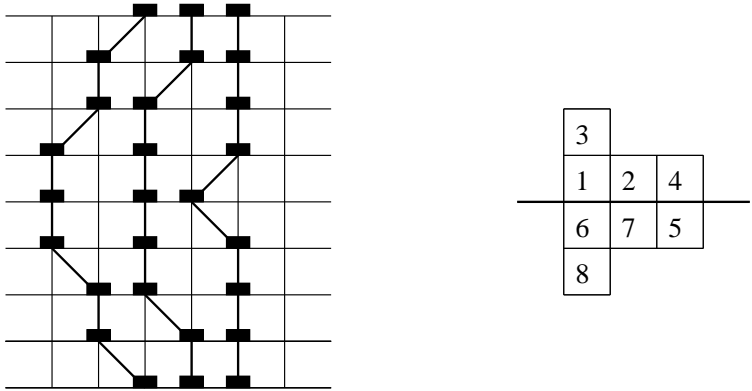}}
\caption{\label{f1} 
\small
A particular configuration of 3 random turns
walkers performing 8 steps in the sequence $L^4 R^4$. The walk can
conveniently be represented diagramatically as done at right
according to the rules specified in the text}
\end{figure}

Our interest is in the random turns vicious walker model \cite{Fi84,Fo91f}.
Here, in the stochastic evolution picture, at each time step $t$
($t = 1,2,\dots$) one particle is selected at random and moved one lattice
site to the left with probability $w_{-1}$, or one lattice site to the right
with probability $w_1$ ($w_{-1} + w_1 = 1$). However, if the lattice site
to the left (right) is already occupied, then the chosen walker moves to the
right (left) with probability one unless this lattice site too is occupied.
In the latter situation another walker is selected at random and the
procedure repeated until one walker has been moved. The move of precisely
one walker so determines the state at time interval $t$. An example of
some typical configurations in the directed paths picture is given in
Figure \ref{f1}. We remark that the random turns vicious walker model can
also be regarded as a particular asymmetric exclusion process \cite{Li85}.

Two aspects of the theory of the random turns vicious walker model are the
subject of this Letter. The first concerns the number of configurations
that have the $p$ walkers initially on adjacent lattice sites
$(l=1,\dots,p)$ on the one-dimensional line, and have the walkers
returned to the same sites after $2n$ time steps (for an odd number of time
steps this is not possible -- thus the reason for $2n$). This will be shown
to be simply related to the number of permutations of $\{1,2,\dots,n\}$ such
that the maximum increasing subsequence has length no greater than $p$.
Then known results concerning the distribution of the maximum increasing
subsequence of a random permuation will be used to determine the
distribution of the maximum displacement of the walkers starting near
the $l=1$ boundary of the initial positions for a variant of the random turns
model. In this variant the number of walkers is greater than or equal to
$n$, and all walkers must move to the left for time steps up to $n$,
while they must move move to the right thereafter, returning to their
starting points at time $2n$.

To count the configurations, label the particles by their initial location
of the one-dimensional line $l=1,\dots,p$. At each time step one walker
will move one lattice site to the left ($L$) or right ($R$), subject to the
rule that no two particles can occupy the same lattice site (in the counting
problem we take $w_1 = w_{-1}$). The constraint that the particles return
to their initial positions after $2n$ steps requires that for each walker
the number of steps to the left equals to the number of steps to the
right after $2n$ steps, and that in total there are $n$ steps $L$ and
$n$ steps $R$. This latter fact allows the walks to be partitioned according
to the ordering of $L$'s and $R$'s, of which there are $\Big (
{2n \atop n} \Big )$ possibilities. The simplest of these is
$L^nR^n$, which means the first $n$ steps are all to the left, while the final
$n$ steps are all to the right. We now pose the question: for a given
combination of $\{L^n, R^n\}$ corresponding to a particular ordering of
$L$, $R$ steps, how many distinct configurations of the $p$ random turns
walkers are there? 

Consider the particular ordering $L^nR^n$. We will represent each
configuration diagramatically by placing a square marked $t$ in column 
$k$ $(k=1,\dots,p)$ if walker $k$ moves at time step $t$. This square must
be placed above any other squares in the column if the step is $L$, and
below any other squares in the column if the step is $R$. If there are no
other squares in the column (i.e.~the walker $k$ is making its first move)
then the square is placed immediately above the axis if the step is
$L$, or immediately below the axis if the step is $R$. An example
of such a diagramatical representation is given in Figure \ref{f1}.

The problem now is to determine the number of distinct diagrams. Consider
the region of the diagram above the axis. Here each column, of which there
are no more than $p$, must be weakly decreasing in height. Each of the 
$n$ squares must be labelled by a different integer $1,\dots,n$ with the
numbers strictly increasing along rows and up columns. We recognize
such a diagram as equivalent to a standard Young tableau (see 
e.g.~\cite{Fu97}; the Young tableau results be interchanging rows and 
columns). Thus there is a bijection between the number of walks
following the sequence $L^n$ and standard Young tableaux with entries
$1,\dots,n$ and no more than $p$ rows. We remark that this is not the first
time a bijection between Young tableaux and vicious walker problems has
been observed: in \cite{GOV98} a bijection between semi-standard tableaux
and the configurations in a special case of the lock-step random walker
model \cite{Fi84, Fo90b} was identified. 

\begin{figure}
\epsfxsize=10cm
\centerline{\epsfbox{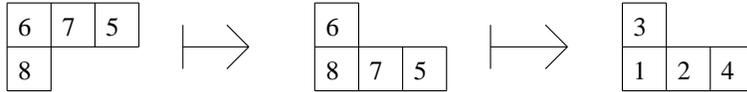}}
\caption{\label{f2} 
\small
Transforming the diagram below the horizontal axis
in Figure \ref{f1} into equivalent Young tableau form}
\end{figure}

Regarding the last $n$ walks, because each walker must return to its starting
position after the $2n$ steps, the shape of the diagram below the $x$-axis
must be the mirror image of the shape of the diagram below the 
$x$-axis. Furthermore, by pushing each column $k=2,\dots,p$ of this
downwards so that they all align with the first entry of the first column,
and then relabelling each box by $t \mapsto 2n + 1 -t$ we see that another
standard Young tableau results. (This procedure is equivalent to constructing
this diagram starting with the final step of the walk and working
backwards.) These transformations are indicated in Figure \ref{f2}.

Thus configurations of the $p$ random turns walkers under consideration
are in bijective correspondence with pairs of standard tableaux of
$n$-boxes constrained so that there are no more than $p$ rows
(recall the role of rows and columns in Young tableaux and our diagrams
is reversed). But such pairs of standard Young tableaux are well known
(see e.g.~\cite{Fu97})
to be in bijective correspondence with permutations of $\{1,\dots,n\}$
such that the length of each increasing subsequence is less than
or equal to $p$. Hence we have enumerated the number of walks in terms
of such permutations.

\begin{prop}\label{prop1}
Consider $p$ random  turns walkers, initially equally spaced one unit
apart and returning to their initial position after $2n$ steps. Suppose
the walkers make their steps in the sequence $L^n R^n$. The total number
of distinct configurations equals the number of permutations of
$\{1,\dots,n\}$ such that the length of the maximum increasing
subsequence is less than or equal to $p$.
\end{prop}

\begin{figure}
\epsfxsize=8cm
\centerline{\epsfbox{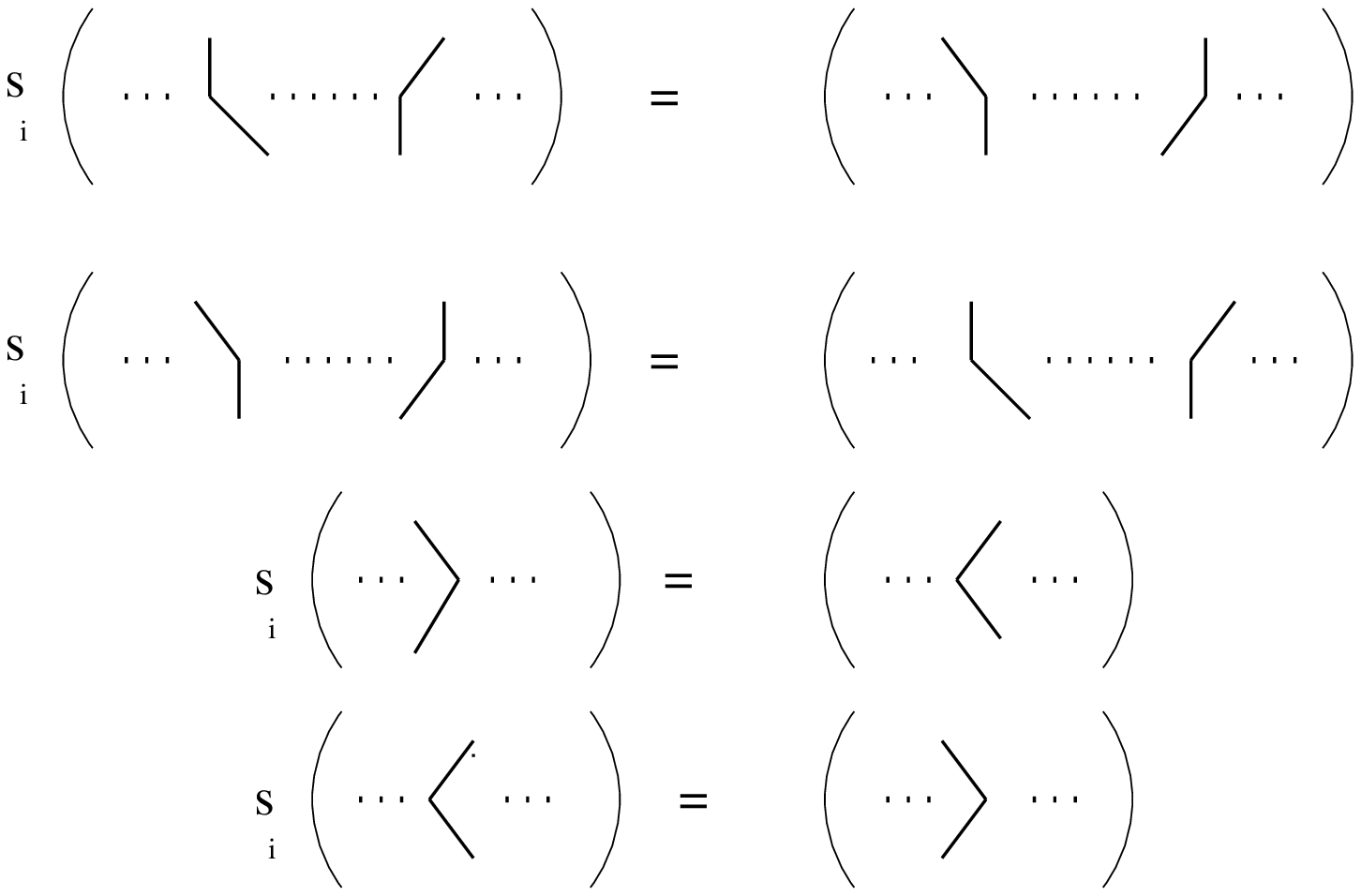}}
\caption{\label{f4} \small
In the first two cases the distinct walkers correspond
to the columns of boxes $i$ and $i+1$, while in the last two cases
boxes $i$ and $i+1$ are in the same column, and the resulting configuration
assumed admissable}
\end{figure}

Consider now another sequence of $n$ $L$'s and $n$ $R$'s. This sequence
can be transformed into the sequence $L^nR^n$ by elementary
transpositions $s_i$ which interchange the $i$th and $(i+1)$th
members of the sequence, assumed to be $R$ and $L$ respectively.
Likewise, we can define the corresponding action of $s_i$ on a diagram
(or equivalently the lattice paths) and so obtain a bijection between
distinct configurations with walks following the sequence $L^nR^n$,
and distinct configurations with walks following
some combination of the sequence $L^nR^n$. We assume step $i$ is opposite
in direction to step $i+1$, and defined the action of $s_i$ to first
interchange the position of boxes $i$ and $i+1$. At the level of the lattice
paths this has the action depicted in Figure \ref{f4}.

\begin{figure}
\epsfxsize=8cm
\centerline{\epsfbox{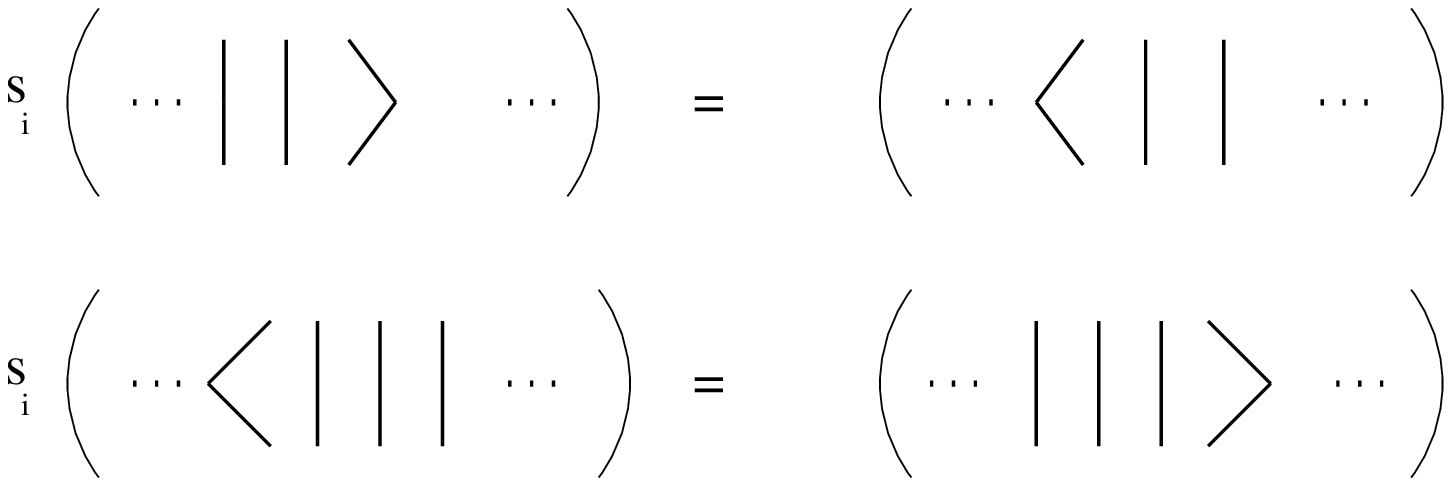}}
\caption{\label{f5} 
\small Here both box $i$ and $i+1$ are in the 
same column and the corresponding walker configuations going from step
$i$ to step $i+1$ are as indicated. Because the rule of Figure \ref{f4}
would lead to inadmissable configurations in each case, $s_i$ acts by
propagating the new configuration to the left and right in the two
cases respectively until an admissable configuration is obtained. }
\end{figure}

If the two boxes are in the same column, it may happen that the resulting
lattice path is inadmissable, in that the new (local) configuration intersects
with an existing path (note that this cannot happen in the first two
cases of Figure \ref{f4}). In such a circumstance we move the
configuration to the left (right) in the second last (last) situation of
Figure \ref{f4} until a permissable configuration is obtained (after so
moving the left-right pairs two vertical lines, corresponding to a
stationary walker, take its place). This is illustrated in Figure
\ref{f5}. Notice that in all cases $s_i^2 = 1$, so the procedure is
invertable. Furthermore, the braid relations $s_is_{i+1}s_i =
s_{i+1}s_is_{i+1}$ are satisfied so the correspondence is
independent of the order of application of the transpositions.

\begin{figure}
\epsfxsize=7cm
\centerline{\epsfbox{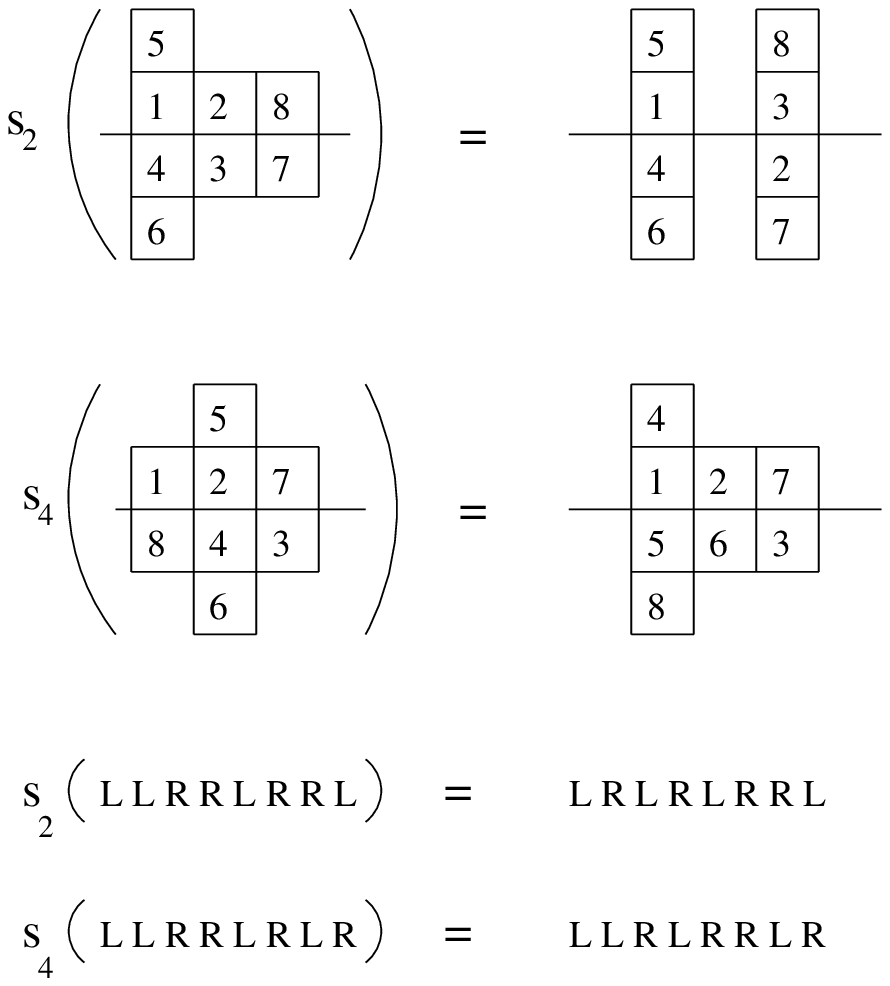}}
\caption{\label{f6} An example of the action of the elementary
transposition operator on some permissable diagrams in the case that
both boxes are from the same column with the situation of 
Figure \ref{f4} applying. The corresponding sequences of $L$'s and
$R$'s are also given. }
\end{figure}

Let's now describe this aspect of the definition of the action of $s_i$
at the level of the diagram. For this purpose we must specify an
admissable diagram. From the rules of the random turns model we see that
an admissable diagram is constructed by putting boxes in the columns
1 to $p$ subject to the conditions that box $i$ must go above (below) the
axis if the step is $L$ ($R$), and that for each column the number of
boxes above the axis minus the number below must be greater than or equal to
the same quantity for the column on the right. Some examples of such
diagrams and the corresponding combinations of $\{L^n, R^n\}$ are given
in Figure \ref{f6}. So we want to describe the situation in which boxes
$i$ and $i+1$ are interchanged but an inadmissable diagram results.
According to Figure \ref{f5} both boxes should be shifted to the left
(right) if the new position of box $i$ is below (above) the axis until
an admissable diagram results. Also we adopt the convention that
numbers in each column should be increasing above and below the axis, so
we rearrange the boxes within a column above and below the axis 
accordingly.

The action of the elementary transpositions thus described give a
bijection between lattice paths with steps in the sequence $L^nR^n$,
and lattice paths with steps in a sequence of any particular
combination of $\{L^n, R^n \}$. Thus by making use of Proposition
\ref{prop1} we can solve the enumeration problem for all such
walks.

\begin{prop}\label{prop2}
The result of Proposition \ref{prop1} for the number of walks in the
sequence $L^n R^n$ also applies for any combination of $\{L^n, R^n\}$.
\end{prop} 

It is of interest to note that there are multiple integral formulas for
both the number of random permutations of $\{1, \dots, n\}$ with at
most $p$ increasing subsequences, and the total number of random turns
paths with $p$ walkers starting at sites $l_1', \dots, l_p'$ and
finishing at sites $l_1,\dots, l_p$ in $2n$ steps. Let us denote these
numbers by $f_{np}$ and $Z_{2n}(l_1', \dots, l_p';
l_1,\dots, l_p)$ respectively. Then we have \cite{Ra98}
\begin{equation}\label{aa}
f_{np} = {(n!)^2 \over (2n)!} {1 \over p!} {1 \over (2 \pi)^p}
\int_{-\pi}^\pi d\theta_1 \cdots \int_{-\pi}^\pi d \theta_p
\Big ( \sum_{j=1}^p 2 \cos \theta_j \Big )^{2n}
\prod_{1 \le \alpha < \beta \le p} |e^{i \theta_\alpha} -
e^{i \theta_\beta} |^2
\end{equation}
and \cite{Fo91f}
\begin{equation}\label{ab}
Z_{2n}(l_1', \dots, l_p';
l_1,\dots, l_p) =
 {1 \over (2 \pi)^p}
\int_{-\pi}^\pi d\theta_1 \cdots \int_{-\pi}^\pi d \theta_p
\Big ( \sum_{j=1}^p 2 \cos \theta_j \Big )^{2n}
\det \Big [ e^{-i (l_\alpha - l_\beta') \theta_\alpha}
\Big ]_{\alpha,\beta = 1,\dots, p}
\end{equation}

Let us specialize (\ref{ab}) to the situation of Propositions
\ref{prop1} and \ref{prop2} by choosing $l_j = l_j' = j$
$(j=1,\dots,p)$. Now 
\begin{eqnarray*}
\det [e^{-i(\alpha - \beta) \theta_\alpha} ]_{\alpha,\beta = 1,\dots, p}
& = & \prod_{j=1}^p e^{-i(j-1) \theta_j} \det [e^{i(\beta - 1) \theta_\alpha}
]_{\alpha,\beta = 1,\dots, p} \\ & = &
\prod_{j=1}^p e^{-i(j-1) \theta_j}
\prod_{1 \le \alpha < \beta \le p} (e^{i\theta_\beta} - e^{i \theta_\alpha}),
\end{eqnarray*}
where the final equality follows from the Vandermonde formula. 
Since the
product of differences is antisymmetric 
the non-symmetric factor in the integrand 
$\prod_{j=1}^p e^{-i (j-1) \theta_j}$ can be anti-symmetrized,
giving another Vandermonde product, provided we divide by  $p!$. Thus we have
\begin{eqnarray}\label{ac}
\lefteqn{
Z_{2n}(\{l_j' = j\}_{j=1,\dots,p};\{l_k = k\}_{k=1,\dots,p})
} \nonumber \\ && 
= {1 \over p!} {1 \over (2 \pi)^p}
\int_{-\pi}^\pi d\theta_1 \cdots \int_{-\pi}^\pi d \theta_p
\Big ( \sum_{j=1}^p 2 \cos \theta_j \Big )^{2n}
\prod_{1 \le \alpha < \beta \le p}
| e^{i\theta_\beta} - e^{i \theta_\alpha} |^2.
\end{eqnarray}
Comparing (\ref{aa}) and (\ref{ac}) gives
\begin{equation}\label{ad}
Z_{2n}(\{l_j' = j\}_{j=1,\dots,p};\{l_k = k\}_{k=1,\dots,p})
= \Big ( {2n \atop n} \Big ) f_{np}
\end{equation}
which is of course also an immediate corollary of Proposition \ref{prop1}
and \ref{prop2}. However once having deduced (\ref{ad}),
the formula (\ref{aa}) for $f_{np}$ follows as a special case of
(\ref{ab}).

We know from the derivation of Proposition \ref{prop1} that there is a
bijection between configurations of $p$ random turns walkers performing
$2n$ steps in the sequence $L^n R^n$ before returning to their initial
positions of all one unit apart, and pairs of standard 
tableaux each of the same shape and
consisting of $n$ boxes. Furthermore, in the bijection the length of
row $j$ corresponds to the maximum displacement of walker $j$ to the left
of its starting point (this occurs at step $n$). Thus if we choose
$p \ge n$ this same bijection holds but with no restriction on the number
of rows. In such a situation the asymptotics of the row lengths are known
precisely \cite{BDJ98,Ok99,BOO99,Jo99b}. 
We can therefore give these results an interpretation
in the random walker setting.

\begin{prop}\label{p3}
Denote the positions of the walkers on the one-dimensional line in the
above asymmetric random turns model by $l_j$ where initially walker
$j$ is at position $l_j = j$. Define the scaled displacements by
\begin{equation}\label{lt}
\tilde{l}_j := n^{1/3}\Big ( {l_j \over n^{1/2}} - 2 \Big )
\end{equation}
and the corresponding scaled $k$-point distribution function by
$$
\rho_k(\tilde{l}_1,\dots,\tilde{l}_k) :=
\lim_{n \to \infty} \Big ( {1 \over n^{1/6}} \Big )^k
\rho_k^{(n)}(\tilde{l}_1,\dots,\tilde{l}_k)
$$
where $\rho_k^{(n)}$ denotes the $k$-point distribution for the walker
problem in the finite system. Then from the results of \cite{Ok99,
BOO99,Jo99b} for the
tableau problem we have
\begin{equation}\label{lt1}
\rho_k(\tilde{l}_1,\dots,\tilde{l}_k) =
\det \Big [ K(\tilde{l}_\alpha, \tilde{l}_\beta) \Big ]_{\alpha, \beta =
1,\dots, k}
\end{equation}
where
$$
K(x,y) := {{\rm Ai}(x) {\rm Ai}'(y) - {\rm Ai}'(x) {\rm Ai}(y) \over x - y}.
$$
\end{prop}

The distribution (\ref{lt1}) is precisely the scaled distribution of the
eigenvalues at the edge of the spectrum for GUE random matrices
\cite{Fo93a,TW94a}. Perhaps more relevantly to the random walker problem,
(\ref{lt}) coincides with the scaled distribution for free fermions on
a line confined by a one-body harmonic potential, at the edge of the support
of the density. The relevance is that there is a well known relationship
between continuous models of non-intersecting walkers and free fermions
(see e.g.~\cite{dN88}).

Regarding some physical features of Proposition \ref{p3}, note from
(\ref{lt}) that the average displacement is $\mu = 2 n^{1/2}$ (for a
recent independent proof of this result see
 \cite{Jo98a}), with standard deviation 
proportional to $(4n)^{1/6} = \mu^\chi$, $\chi = 1/3$. As emphasized
in \cite{Jo99a}, the exponent $\chi = 1/3$ is typical of two-dimensional
growth models
(it can be derived from the one-dimensional Burgers equation
describing such processes \cite{vB91}). 
On this point we recall that vicious walker paths
fixed at the endpoints as in Proposition \ref{p3} form the well-known
(see e.g.~\cite{JEB91}) terrace-step-kink model of a crystal surface.

\subsection*{Acknowledgement}
The financial support of the ARC, including funds to support the visit
of G.~Olshanski whose lectures benefitted the present work, are
acknowledged. Also, the remarks of T.H.~Baker on the original manuscript
are appreciated. 


\end{document}